\documentclass[12pt, a4paper]{article}
\usepackage{amsmath}
\usepackage{amssymb}
\usepackage{t1enc}
\usepackage[latin1]{inputenc}
\usepackage[english]{babel}
\usepackage{amsfonts}
\usepackage{latexsym}
\usepackage{bm}
\newtheorem{theorem}{Theorem}[section]
\newtheorem{prop}[theorem]{Proposition}

\newtheorem{cor}[theorem]{Corollary}
\newtheorem{conj}[theorem]{Conjecture}

\begin{document}

\title{Intersecting families of sets and permutations: a survey}

\author{Peter Borg\\[5mm]
Department of Mathematics, University of Malta, \\
Msida MSD 2080, Malta\\
\texttt{p.borg.02@cantab.net}}
\date{\today} \maketitle

\begin{abstract}
A family $\mathcal{A}$ of sets is said to be
\emph{$t$-intersecting} if any two sets in $\mathcal{A}$ have at
least $t$ common elements. A central problem in extremal set
theory is to determine the size or structure of a largest
$t$-intersecting sub-family of a given family $\mathcal{F}$. We
give a survey of known results, conjectures and open problems for
various important families $\mathcal{F}$, namely, power sets,
levels of power sets, hereditary families, families of signed
sets, families of labeled sets, and families of permutations. We
also provide some extensions and consequences of known results.
\end{abstract}

\section{Introduction}

%
Unless otherwise stated, we shall use small letters such as $x$ to
denote elements of a set or non-negative integers or functions,
capital letters such as $X$ to denote sets, and calligraphic
letters such as $\mathcal{F}$ to denote \emph{families}
(i.e.~sets whose elements are sets themselves). It is to be
assumed that arbitrary sets and families are \emph{finite}. We
call a set $A$ an \emph{$r$-element set}, or simply an
\emph{$r$-set}, if its size $|A|$ is $r$ (i.e.~if it contains
exactly $r$ elements). A family is said to be \emph{uniform} if
all its sets are of the same size. 

The set $\{1, 2, ...\}$ of positive integers is denoted by
$\mathbb{N}$. For $m, n \in \mathbb{N}$ with $m \leq n$, the set
$\{i \in \mathbb{N} \colon m \leq i \leq n\}$ is denoted by
$[m,n]$, and if $m = 1$ then we also write $[n]$. For a set $X$,
the \emph{power set} $\{A \colon A \subseteq X\}$ of $X$ is
denoted by $2^X$, and the uniform sub-family $\{Y \subseteq
X\colon |Y| = r\}$ of $2^X$ is denoted by $X \choose r$.

For a family $\mathcal{F}$ of sets, we denote the union of all
sets in $\mathcal{F}$ by $U(\mathcal{F})$ and we denote the size
of a largest set in $\mathcal{F}$ by $\alpha(\mathcal{F})$. For an
integer $r \geq 0$, we denote the uniform sub-family $\{F \in
\mathcal{F} \colon |F| = r\}$ of $\mathcal{F}$ by
$\mathcal{F}^{(r)}$ (note that $\mathcal{F}^{(r)} = {X \choose
r}$ if $\mathcal{F} = 2^X$), and we call $\mathcal{F}^{(r)}$ the
\emph{$r$'th level of $\mathcal{F}$}. For a set $S$, we denote
$\{F \in \mathcal{F} \colon S \subseteq F\}$ by $\mathcal{F}(S)$.
We may abbreviate $\mathcal{F}(\{x\})$ to $\mathcal{F}(x)$. If $x
\in U(\mathcal{F})$ then we call $\mathcal{F}(x)$ a \emph{star of
$\mathcal{F}$}. More generally, if $T$ is a $t$-element subset of
a set in $\mathcal{F}$, then we call $\mathcal{F}(T)$ a
\emph{$t$-star of $\mathcal{F}$}.

A family $\mathcal{A}$ is said to be \emph{intersecting} if $A
\cap B \neq \emptyset$ for any $A, B \in \mathcal{A}$.  More
generally, $\mathcal{A}$ is said to be \emph{$t$-intersecting} if
$|A \cap B| \geq t$ for any $A, B \in \mathcal{A}$. So an
intersecting family is a $1$-intersecting family. A
$t$-intersecting family $\mathcal{A}$ is said to be
\emph{trivial} if $|\bigcap_{A \in \mathcal{A}}A| \geq t$
(i.e.~there are at least $t$ elements common to all the sets in
$\mathcal{A}$); otherwise, $\mathcal{A}$ is said to be
\emph{non-trivial}. So a $t$-star of a family $\mathcal{F}$ is a
trivial $t$-intersecting sub-family of $\mathcal{F}$ that is not
contained in any other. If there exists a $t$-set $T$ such that
$\mathcal{F}(T)$ is a largest $t$-intersecting sub-family of
$\mathcal{F}$ (i.e.~no $t$-intersecting sub-family of
$\mathcal{F}$ has more sets than $\mathcal{F}(T)$), then we say
that $\mathcal{F}$ has the \emph{$t$-star property at $T$}, or we
simply say that $\mathcal{F}$ has the \emph{$t$-star property}.
If either $\mathcal{F}$ has no $t$-intersecting sub-families
(which is the case if and only if $\alpha(\mathcal{F}) < t$) or
all the largest $t$-intersecting sub-families of $\mathcal{F}$ are
$t$-stars, then we say that $\mathcal{F}$ has the \emph{strict}
$t$-star property. We may abbreviate `$1$-star property' to `star
property'.


\emph{Extremal set theory} is the study of how small or how large
a system of sets can be under certain conditions. In this paper
we are concerned with the following natural and central problem
in this field.\\
\\
\textbf{Problem:} \emph{Given a family $\mathcal{F}$ and an
integer $t \geq 1$, determine the size or structure of a largest
$t$-intersecting sub-family of $\mathcal{F}$.}\\
\\
We provide a survey of results that answer this question for
families that are of particular importance, and we also point out
open problems and conjectures. The survey papers \cite{DF} and
\cite{F} cover a few of the results we mention here and also go
into many variations of the above problem, however, much progress
has been made since their publication. Here we cover many of the
important results that have been established to date, restricting
ourselves to the problem stated above.

The most obvious families to consider are the power set $2^{[n]}$
and the uniform sub-family ${[n] \choose r}$, and in fact the
problem for these families has been solved completely. However,
there are other important families on which much progress has
been made, and there are others that are still subject to much
investigation. The families defined below are perhaps the ones
that have received most attention and that we will be concerned
with.\\
\\
\textbf{Hereditary families:} A family $\mathcal{H}$ is said to be
a \emph{hereditary family} (also called an \emph{ideal} or a
\emph{downset}) if all the subsets of any set in $\mathcal{H}$
are in $\mathcal{H}$. Clearly a family is hereditary if and only
if it is a union of power sets. A \emph{base} of $\mathcal{H}$ is
a set in $\mathcal{H}$ that is not a subset of any other set in
$\mathcal{H}$. So a hereditary family is the union of power sets
of its bases. An example of a hereditary family is the family of
independent sets of a \emph{graph} or \emph{matroid}.\\
\\
\textbf{Families of signed sets:} Let $X$ be an $r$-set $\{x_1,
..., x_r\} \subset \mathbb{N}$. Let $y_1, ..., y_r \in
\mathbb{N}$. We call the set $\{(x_1,y_1), ..., (x_r,y_r)\}$ a
\emph{$k$-signed $r$-set} if $\max\{y_i \colon i \in [r]\} \leq
k$. For an integer $k \geq 2$ we define $\mathcal{S}_{X,k}$ to be
the family of $k$-signed $r$-sets given by
\[\mathcal{S}_{X,k} := \{\{(x_1,y_1), ..., (x_r,y_r)\} \colon y_1,
..., y_r \in [k]\}.\]
So a set $A$ is a member of $\mathcal{S}_{X,k}$ if and only if it is a subset of the Cartesian product $X \times [k] := \{(x,y) \colon x
\in X, y \in [k]\}$ satisfying $|A \cap (\{x\} \times [k])| = 1$ for all $x \in X$. We shall set $\mathcal{S}_{\emptyset,k} := \emptyset$. With a slight abuse of notation, for a family $\mathcal{F}$ we define
\[\mathcal{S}_{\mathcal{F},k} := \bigcup_{F \in
\mathcal{F}}\mathcal{S}_{F,k}.\]
\\
\textbf{Families of labeled sets:} For ${\bf k} := (k_1, ...,
k_n)$ with $k_1, ..., k_n \in \mathbb{N}$ and $k_1 \leq ... \leq
k_n$, we define the family $\mathcal{L}_{\bf{k}}$ of \emph{labeled
$n$-sets} by
\[\mathcal{L}_{\bf{k}} := \{\{(1,y_1), ..., (n,y_n)\} \colon y_i
\in [k_i] \mbox{ for each } i \in [n]\}.\]
Note that $\mathcal{S}_{[n],k} = \mathcal{L}_{(k_1, ..., k_n)}$
with $k_1 = ... = k_n = k$.

An equivalent formulation for $\mathcal{L}_{\bf{k}}$ is the
Cartesian product $[k_1] \times ... \times [k_n] := \{(y_1, ...,
y_n) \colon y_i \in [k_i] \mbox{ for each } i \in [n]\}$, but it
is more convenient to work with $n$-sets than work with
$n$-tuples (the alternative formulation demands that we change
the setting of families of sets to one of sets of $n$-tuples).

For any $r \in [n]$, we define
\[\mathcal{L}_{{\bf k},r} := \{\{(x_1,y_{x_1}), ...,
(x_r,y_{x_r})\} \colon \{x_1, ..., x_r\} \in {[n] \choose r}, \,
y_{x_i} \in [k_{x_i}] \mbox{ for each } i \in [r]\},\]
and we set $\mathcal{L}_{{\bf k},0} = \emptyset$. Thus, for any
$0 \leq r \leq n$, $\mathcal{L}_{{\bf k},r}$ is the family of all
$r$-element subsets of all the sets in $\mathcal{L}_{\bf{k}}$, and
$\mathcal{L}_{{\bf k},n} = \mathcal{L}_{\bf{k}}$. We also define
$\mathcal{L}_{{\bf k}, \leq r} := \bigcup_{i=0}^r
\mathcal{L}_{{\bf k},i}$.\\
\\
\textbf{Families of permutations:} For an $r$-set $X := \{x_1,
..., x_r\}$, we define $\mathcal{S}_{X,k}^*$ to be the special
sub-family of $\mathcal{S}_{X,k}$ given by
\[\mathcal{S}_{X,k}^* := \left\{\{(x_1, y_1), ..., (x_r, y_r)\}
\colon y_1, ..., y_r \mbox{ are \emph{distinct} elements of }
[k]\right\}.\]
Note that $\mathcal{S}_{X,k}^* \neq \emptyset$ if and only if $r
\leq k$. With a slight abuse of notation, for a family
$\mathcal{F}$ we define $\mathcal{S}_{\mathcal{F},k}^*$ to be the
special sub-family of $\mathcal{S}_{\mathcal{F},k}$ given by
$$\mathcal{S}_{\mathcal{F},k}^* := \bigcup_{F \in
\mathcal{F}}\mathcal{S}_{F,k}^*.$$

An \emph{r-partial permutation of a set $N$} is a pair $(A, f)$
where $A \in {N \choose r}$ and $f \colon A \rightarrow N$ is an
injection. An $|N|$-partial permutation of $N$ is simply called a
\emph{permutation of $N$}. Clearly, the family of permutations of
$[n]$ can be re-formulated as $\mathcal{S}_{[n],n}^*$, and the
family of $r$-partial permutations of $[n]$ can be re-formulated
as $\mathcal{S}_{{[n] \choose r},n}^*$.

Let $X$ be as above. $\mathcal{S}_{X,k}^*$ can be interpreted as
the family of permutations of sets in ${[k] \choose r}$: consider
the bijection $\beta \colon \mathcal{S}_{X,k}^* \rightarrow
\{(A,f) \colon A \in {[k] \choose r}, f \colon A \rightarrow A
\mbox{ is a bijection}\}$ defined by $\beta(\{(x_1,a_1), ...,$
$(x_r,a_r)\}) := (\{a_1, ..., a_r\}, f)$ where, for $b_1 < ... <
b_r$ such that $\{b_1, ..., b_r\} = \{a_1, ..., a_r\}$, $f(b_i) :=
a_i$ for $i = 1, ..., r$. $\mathcal{S}_{X,k}^*$ can also be
interpreted as the sub-family $\mathcal{X} := \{(A,f) \colon A
\in {[k] \choose r}, f \colon A \rightarrow [r] \mbox{ is a
bijection}\}$ of the family of $r$-partial permutations of $[k]$:
consider an obvious bijection from $\mathcal{S}_{X,k}^*$ to
$\mathcal{S}_{{[k] \choose r},r}^*$ and another one from
$\mathcal{S}_{{[k] \choose r},r}^*$ to $\mathcal{X}$.


\section{Intersecting sub-families of ${[n] \choose r}$ and
$2^{[n]}$} \label{powerset}

The study of intersecting families took off with the publication
of \cite{EKR}, which features the following classical result,
known as the Erd\H os-Ko-Rado (EKR) Theorem.

\begin{theorem}[EKR Theorem \cite{EKR}] \label{EKR} If $r \leq
n/2$ and $\mathcal{A}$ is an intersecting sub-family of ${[n]
\choose r}$, then $|\mathcal{A}| \leq {n-1 \choose r-1}$.
\end{theorem}
This means that for $r \leq n/2$, ${[n] \choose r}$ has the star
property, because the bound ${n-1 \choose r-1}$ is the size of any
star of ${[n] \choose r}$. Note that if $r > n/2$, then any two
$r$-element subsets of $[n]$ must intersect, and hence ${[n]
\choose r}$ is an intersecting family (also note it is a
non-trivial one, so ${[n] \choose r}$ does not have the star
property in this case).

In order to prove Theorem~\ref{EKR}, Erd\H os, Ko and Rado
\cite{EKR} introduced a method known as \emph{compression} or
\emph{shifting}; see \cite{F} for a survey on the uses of this
powerful technique in extremal set theory. There are various
proofs of Theorem~\ref{EKR}, two of which are particularly short
and beautiful: Katona's \cite{K} proof, which featured an elegant
argument known as the \emph{cycle method}, and Daykin's proof
\cite{D} using another fundamental result known as the
Kruskal-Katona Theorem \cite{Ka,Kr}. Hilton and Milner \cite{HM}
proved that for $r \leq n/2$, the family $\mathcal{N}_{n,r} :=
\left\{A \in {[n] \choose r} \colon A \cap [2,r+1] \neq
\emptyset\right\} \cup \{[2,r+1]\}$ is a largest non-trivial
intersecting sub-family of ${[n] \choose r}$, and since the size
of $\mathcal{N}_{n,r}$ is ${n-1 \choose r-1} - {n-r-1 \choose r-1}
+ 1$, it follows that if $r < n/2$, then the stars of ${[n]
\choose r}$ are the largest intersecting sub-families of ${[n]
\choose r}$, i.e.~${[n] \choose r}$ has the strict star property.
Note that if $r = n/2$, then any sub-family $\mathcal{A}$ of
${[n] \choose r}$ satisfying $|\mathcal{A} \cap \{A, [2r]
\backslash A\}| = 1$ for all $A \in {[n] \choose r}$ is an
intersecting sub-family of ${[n] \choose r}$ of size
$\frac{1}{2}{n \choose r} = \frac{1}{2}{2r \choose r} = {2r-1
\choose r-1}$, and hence one of maximum size (an example of such a
family $\mathcal{A}$ is $\mathcal{N}_{2r,r}$, so ${[n] \choose r}$
does not have the strict star property if $r = n/2$).

Also in \cite{EKR}, Erd\H os, Ko and Rado initiated the study of
$t$-intersecting families for $t \geq 2$. They proved that for $t
< r$, there exists an integer $n_0(r,t)$ such that for all $n \geq
n_0(r,t)$, the largest $t$-intersecting sub-families of ${[n]
\choose r}$ are the $t$-stars (which are of size ${n-t \choose
r-t}$). For $t \geq 15$, Frankl \cite{F_t1} showed that the
smallest such $n_0(r,t)$ is $(r-t+1)(t+1) + 1$ and that if $n =
(r-t+1)(t+1)$, then ${[n] \choose r}$ still has the $t$-star
property but not the strict $t$-star property. Subsequently,
using algebraic means, Wilson \cite{W} proved that ${[n] \choose
r}$ has the $t$-star property for any $t$ and $n \geq
(r-t+1)(t+1)$. Frankl \cite{F_t1} conjectured that among the
largest $t$-intersecting sub-families of ${[n] \choose r}$ there
is always at least one of the families $\left\{A \in {[n] \choose
r} \colon |A \cap [t+2i]| \geq t+i \right\}$, $i = 0, 1, ...,
r-t$. A remarkable proof of this long-standing conjecture
together with a complete characterisation of the extremal
structures was finally obtained by Ahlswede and Khachatrian
\cite{AK1} by means of the \emph{compression} technique introduced
in \cite{EKR}.

\begin{theorem}[Ahlswede and Khachatrian \cite{AK1}]\label{AK} Let
$1 \leq t \leq r \leq n$ and let $\mathcal{A}$ be a largest
$t$-intersecting sub-family of ${[n] \choose r}$. \\
(i) If $(r-t+1)(2 + \frac{t-1}{i+1}) < n <
(r-t+1)(2+\frac{t-1}{i})$ for some $i \in \{0\} \cup \mathbb{N}$ -
where, by convention, $(t-1)/i = \infty$ if $i = 0$ - then
$\mathcal{A} = \{A \in {[n] \choose r} \colon |A
\cap X| \geq t+i\}$ for some $X \in {[n] \choose t+2i}$. \\
(ii) If $t \geq 2$ and $(r-t+1)(2 + \frac{t-1}{i+1}) = n$ for some
$i \in \{0\} \cup \mathbb{N}$, then $\mathcal{A} = \{A \in {[n]
\choose r} \colon |A \cap X| \geq t+j\}$ for some $j \in
\{i,i+1\}$ and $X \in {[n] \choose t+2j}$.
\end{theorem}
It is worth mentioning that in \cite{AK3} Ahlswede and Khachatrian
went on to determine the largest non-trivial $t$-intersecting
sub-families of ${[n] \choose r}$.

Erd\H os, Ko and Rado \cite{EKR} pointed out the simple fact that
$2^{[n]}$ has the star property (indeed, for any set $A$ in an
intersecting sub-family $\mathcal{A}$ of $2^{[n]}$, the
complement $[n] \backslash A$ cannot be in $\mathcal{A}$, and
hence the size of $\mathcal{A}$ is at most $\frac{1}{2}|2^{[n]}|
= 2^{n-1}$, i.e.~the size of a star of $2^{[n]}$); note that there
are many non-trivial intersecting sub-families of $2^{[n]}$ of
maximum size $2^{n-1}$ (such as $\{A \subseteq [n] \colon |A \cap
[3]| \geq 2\}$), so $2^{[n]}$ does not have the strict star
property. They also asked what the size of a largest
$t$-intersecting sub-family of $2^{[n]}$ is for $t \geq 2$. The
answer in a complete form was given by Katona \cite{Kat}.

\begin{theorem}[Katona \cite{Kat}]\label{K_thm} Let $t \geq 2$,
and let $\mathcal{A}$ be a largest $t$-intersecting sub-family of
$2^{[n]}$. \\
(i) If $n + t = 2l$ then $\mathcal{A} = \{A \subseteq [n] \colon
|A| \geq l\}$. \\
(ii) If $n + t = 2l+1$ then $\mathcal{A} = \{A \subseteq [n]
\colon |A \cap X| \geq l\}$ for some $X \in {[n] \choose n-1}$.
\end{theorem}
It is interesting that for $n > t \geq 2$, $2^{[n]}$ does not
have the $t$-star property.

Many other beautiful results were inspired by the seminal paper
\cite{EKR}, as are the results we present in the subsequent
sections.

\section{Intersecting sub-families of hereditary families}
\label{Hersection}

Recall that $2^{[n]}$ has the star property. Also recall that the
power set of a set $X$ is the simplest example of a hereditary
family as $2^X$ is a hereditary family with only one base ($X$).
An outstanding open problem in extremal set theory is the
following conjecture (see \cite{Borg3} for a more general
conjecture).

\begin{conj}[\cite{Chv}] \label{Chvatal} If $\mathcal{H}$
is a hereditary family, then $\mathcal{H}$ has the star property.
\end{conj}
Chv\'{a}tal \cite{Chva} verified this conjecture for the case
when $\mathcal{H}$ is \emph{left-compressed} (i.e.~$\mathcal{H}
\subseteq 2^{[n]}$ and $(H \backslash \{j\}) \cup \{i\} \in
\mathcal{H}$ whenever $1 \leq i < j \in H \in \mathcal{H}$ and $i
\notin H$). Snevily \cite{Sn} took this result (together with
results in \cite{Schonheim, Wang}) a significant step forward by
verifying Conjecture~\ref{Chvatal} for the case when
$\mathcal{H}$ is \emph{compressed with respect to an element $x$
of $U(\mathcal{H})$} (i.e.~$(H \backslash \{h\}) \cup \{x\} \in
\mathcal{H}$ whenever $h \in H \in \mathcal{H}$ and $x \notin H$).

\begin{theorem} [\cite{Sn}] \label{Snevily} If a hereditary family
$\mathcal{H}$ is compressed with respect to an element $x$ of
$U(\mathcal{H})$, then $\mathcal{H}$ has the star property at $\{x\}$.
\end{theorem}
%
A generalisation is proved in \cite{Borg3} by means of an
alternative self-contained argument. Snevily's proof of
Theorem~\ref{Snevily} makes use of the following interesting
result of Berge \cite{Berge} (a proof of which is also provided
in \cite[Chapter 6]{Anderson}).
\begin{theorem} [\cite{Berge}] \label{Berge} If
$\mathcal{H}$ is a hereditary family, then $\mathcal{H}$ is a
disjoint union of pairs of disjoint sets, together with
$\emptyset$ if $|\mathcal{H}|$ is odd.
\end{theorem}
This result was also motivated by Conjecture~\ref{Chvatal} as it
has the following immediate consequence.

\begin{cor}\label{Bergecor} If $\mathcal{A}$ is an intersecting
sub-family of a hereditary family $\mathcal{H}$, then
$$|\mathcal{A}| \leq \frac{1}{2}|\mathcal{H}|.$$
\end{cor}
\textbf{Proof.} For any pair of disjoint sets, at most only one
set can be in an intersecting family $\mathcal{A}$. By
Theorem~\ref{Berge}, the result follows.~\hfill{$\Box$} \\

A special case of Theorem~\ref{Snevily} is a result of Schönheim
\cite{Schonheim} which says that Conjecture~\ref{Chvatal} is true
if the bases of $\mathcal{H}$ have a common element, and this
follows immediately from Corollary~\ref{Bergecor} and the
following fact.

\begin{prop} [\cite{Schonheim}] \label{Bergeprop} If the bases of
a hereditary family $\mathcal{H}$ have a common element $x$, then
$$|\mathcal{H}(x)| = \frac{1}{2}|\mathcal{H}|.$$
\end{prop}
\textbf{Proof.} Partition $\mathcal{H}$ into $\mathcal{A} :=
\mathcal{H}(x)$ and $\mathcal{B} := \{B \in \mathcal{H} \colon x
\notin B\}$. If $A \in \mathcal{A}$ then $A \backslash \{x\} \in
\mathcal{B}$; so $|\mathcal{A}| \leq |\mathcal{B}|$. If $B \in
\mathcal{B}$ then $B \subseteq C$ for some base $C$ of
$\mathcal{H}$, and hence $B \cup \{x\} \in \mathcal{A}$ since $x
\in C$; so $|\mathcal{B}| \leq |\mathcal{A}|$. Thus
$|\mathcal{A}| = |\mathcal{B}| =
\frac{1}{2}|\mathcal{H}|$.~\hfill{$\Box$}\\

Many other results and problems have been inspired by
Conjecture~\ref{Chvatal} or are related to it; see
\cite{Chvatalsite, Miklos, West}.

Conjecture~\ref{Chvatal} cannot be generalised to the
$t$-intersection case. Indeed, if $n > t \geq 2$ and $\mathcal{H}
= 2^{[n]}$, then by Theorem~\ref{K_thm}, $\mathcal{H}$ does not
have the $t$-star property.

We now turn our attention to uniform intersecting sub-families of
hereditary families, or rather intersecting sub-families of
levels of hereditary families. For any hereditary family
$\mathcal{H}$, let $\mu(\mathcal{H})$ denote the size of a
smallest base of $\mathcal{H}$.

A \emph{graph} $G$ is a pair $(V,E)$ with $E \subseteq {V \choose
2}$, and a set $I \subseteq V$ is said to be an \emph{independent
set of $G$} if $\{i,j\} \notin E$ for any $i, j \in I$. Let
$\mathcal{I}_G$ denote the family of all independent sets of a
graph $G$. Clearly $\mathcal{I}_G$ is a hereditary family.
Holroyd and Talbot \cite{HT} made a nice conjecture which claims
that if $G$ is a graph and $\mu(\mathcal{I}_G) \geq 2r$, then
${\mathcal{I}_G}^{(r)}$ has the star property, and
${\mathcal{I}_G}^{(r)}$ has the strict star property if
$\mu(\mathcal{I}_G) > 2r$. In \cite{Borg} the author conjectured
that this is true for any hereditary family and that in general
the following holds.
\begin{conj}[\cite{Borg}] \label{AK gen} If $t \leq r$,
$\emptyset \neq S \subseteq [t,r]$ and $\mathcal{H}$ is a
hereditary family with $\mu(\mathcal{H}) \geq (t+1)(r-t+1)$,
then: \\
(i) $\bigcup_{s \in S}\mathcal{H}^{(s)}$ has the $t$-star
property;\\
(ii) $\bigcup_{s \in S}\mathcal{H}^{(s)}$ has the strict $t$-star
property if either $\mu(\mathcal{H}) > (t+1)(r-t+1)$ or $S \neq
\{r\}$.
\end{conj}
Note that Theorem~\ref{AK} solves the special case when
$\mathcal{H} = 2^{[n]}$ and tells us that we cannot improve the
condition $\mu(\mathcal{H}) \geq (t+1)(r-t+1)$. The author
\cite{Borg} proved that this conjecture is true if
$\mu(\mathcal{H})$ is sufficiently large.

\begin{theorem}[\cite{Borg}] \label{t int her} Conjecture~\ref{AK gen}
is true if $\mu(\mathcal{H}) \geq (r-t){3r-2t-1 \choose t+1} + r$.
\end{theorem}
The motivation behind establishing this result for any union of
levels of a hereditary family $\mathcal{H}$ within a certain
range is that this general form cannot be immediately deduced
from the result for just one level of $\mathcal{H}$ (i.e.~the case
$S = \{r\}$ in Conjecture~\ref{AK gen}). As demonstrated in
Example~1 in \cite{Borg}, the reason is simply that if $T$ is a
$t$-set such that $\mathcal{H}^{(s)}(T)$ ($s \in [t,r]$) is a
largest $t$-star of the level $\mathcal{H}^{(s)}$, then for $p
\neq s$ ($p \in [t,r]$), $\mathcal{H}^{(p)}(T)$ not only may not
be a largest $t$-star of the level $\mathcal{H}^{(p)}$ but may be
smaller than some non-trivial $t$-intersecting sub-family of
$\mathcal{H}^{(p)}$. This is in fact one of the central
difficulties arising from any EKR-type problem for hereditary
families. In the proof of Theorem~\ref{t int her}, this obstacle
was overcome by showing that for any non-trivial $t$-intersecting
sub-family $\mathcal{A}$ of the union, we can construct a
$t$-star that is larger than $\mathcal{A}$ (and that is not
necessarily a largest $t$-star). Many other proofs of EKR-type
results are based on determining at least one largest $t$-star;
as in the case of each theorem mentioned in
Section~\ref{powerset}, the setting is often symmetrical to the
extent that all $t$-stars are of the same size and of a known
size.

An interesting immediate consequence of Theorem~\ref{t int her} is
that the union of the first $r \geq t$ levels of a hereditary
family $\mathcal{H}$ has the strict $t$-star property if
$\mu(\mathcal{H})$ is sufficiently larger than $r$.

\begin{cor}[\cite{Borg}] If $t \leq r$ and $\mathcal{H}$ is a
hereditary family with $\mu(\mathcal{H}) \geq (r-t){3r-2t-1
\choose t+1} + r$, then $\bigcup_{s=0}^r \mathcal{H}^{(s)}$ has
the strict $t$-star property.
\end{cor}
\textbf{Proof.} Let $\mathcal{A}$ be a $t$-intersecting sub-family
of $\bigcup_{s=0}^r \mathcal{H}^{(s)}$. Then no set in
$\mathcal{A}$ is of size less than $t$, so $\mathcal{A} \subseteq
\bigcup_{s \in S} \mathcal{H}^{(s)}$ with $S = [t,r]$. The result
follows by Theorem~\ref{t int her}.~\hfill{$\Box$}\\

This means that for the special case $t=1$, we have the following.

\begin{cor}[\cite{Borg}] Conjecture~\ref{Chvatal} is true if
$\mathcal{H} = \bigcup_{s=0}^r \mathcal{J}^{(s)}$ for some $r \in
\mathbb{N}$ and some hereditary family $\mathcal{J}$ with
$\mu(\mathcal{J}) \geq \frac{3}{2}(r-1)^2(3r-4) + r$.
\end{cor}
%

The following extension of Theorem~\ref{AK} was also proved in
\cite{Borg1}.

\begin{theorem}[\cite{Borg}] \label{AKext} Conjecture~\ref{AK gen} is true
if $\mathcal{H}$ is left-compressed.
\end{theorem}

\section{Intersecting families of signed sets} \label{SS_Int}

The `signed sets' terminology was introduced in \cite{BL} for a
setting that can be re-formulated as $\mathcal{S}_{{[n] \choose
r},k}$, and the general formulation $\mathcal{S}_{\mathcal{F},k}$
was introduced in \cite{Borg1}, the theme of which is the
following conjecture.
\begin{conj}[\cite{Borg1}]\label{ss conj} For any family
$\mathcal{F}$ and any $k \geq 2$, \\
(i) $\mathcal{S}_{\mathcal{F},k}$ has the star property;\\
(ii) $\mathcal{S}_{\mathcal{F},k}$ does not have the strict
star property only if $k = 2$ and there exist at least
three elements $u_1, u_2, u_3$ of $U(\mathcal{F})$ such that
$\mathcal{F}(u_1) = \mathcal{F}(u_2) = \mathcal{F}(u_3)$ and
$\mathcal{S}_{\mathcal{F},2}((u_1,1))$ is a largest star of
$\mathcal{S}_{\mathcal{F},2}$.
\end{conj}
The converse of (ii) is true, and the proof is simply that $\{A
\in \mathcal{S}_{\mathcal{F},2} \colon |A \cap \{(u_1,1), (u_2,1),
(u_3,1)\}| \geq 2\}$ is a non-trivial intersecting sub-family of
$\mathcal{S}_{\mathcal{F},2}$ that is as large as
$\mathcal{S}_{\mathcal{F},2}((u_1,1))$.

In \cite{Borg3} a similarity between the intersection problem for
hereditary families and the one presented above is demonstrated,
and in fact a conjecture generalising both
Conjecture~\ref{Chvatal} and the above conjecture is suggested.

Recall that a family $\mathcal{F}$ is said to be compressed with
respect to an element $x$ of $U(\mathcal{F})$ if $(F \backslash
\{u\}) \cup \{x\} \in \mathcal{F}$ whenever $u \in F \in
\mathcal{F}$ and $x \notin F$. The following is the main result
in the paper featuring the above conjecture.

\begin{theorem}[\cite{Borg1}] Conjecture~\ref{ss conj} is true if
$\mathcal{F}$ is compressed with respect to an element $x$ of
$U(\mathcal{F})$, and $\mathcal{S}_{\mathcal{F},k}$ has the star
property at $\{(x,1)\}$.
\end{theorem}
Since ${[n] \choose r}$ is compressed with respect to any element
of $[n]$, the above result has the following immediate
consequence, which is a well-known result that was first stated
by Meyer \cite{Meyer} and proved in different ways by Deza and
Frankl \cite{DF}, Bollob\'{a}s and Leader \cite{BL}, Engel
\cite{E} and Erd\H os et al. \cite{EFK}.

\begin{theorem}[\cite{BL,DF,E,EFK}]\label{DFBL} Let $r \in [n]$
and let $k \geq 2$. Then: \\
(i) $\mathcal{S}_{{[n] \choose r},k}$ has the star property;\\
(ii) if $(r,k) \neq (n,2)$ then $\mathcal{S}_{{[n] \choose r},k}$
has the strict star property.
\end{theorem}
Thus the size of an intersecting sub-family of $\mathcal{S}_{{[n]
\choose r},k}$ is at most ${n-1 \choose r-1}k^{r-1}$, i.e.~the
size of any star of $\mathcal{S}_{{[n] \choose r},k}$. Berge
\cite{B} and Livingston \cite{L} had proved (i) and (ii),
respectively, for the special case $\mathcal{F} = \{[n]\}$ (other
proofs are found in \cite{G,M}).

In \cite{Borg1} Conjecture~\ref{ss conj} is also verified for the
case when $\mathcal{F}$ is uniform and has the star property;
Holroyd and Talbot \cite{HT} had essentially proved part (i) of
the conjecture for such a family $\mathcal{F}$ in a
graph-theoretical context.

The $t$-intersection problem for sub-families of
$\mathcal{S}_{[n],k}$ has also been solved. Frankl and Füredi
were the first to investigate it. In \cite{FF2} they conjectured
that among the largest $t$-intersecting sub-families of
$\mathcal{S}_{[n],k}$ there is always one of the families
$\mathcal{A}_i := \{A \in \mathcal{S}_{[n],k} \colon |A \cap ([t +
2i] \times [1])| \geq t + i\}$, $i = 0, 1, 2, ...$, and they
proved that if $k \geq t+1 \geq 16$, then $\mathcal{A}_0$ is
extremal and hence $\mathcal{S}_{[n],k}$ has the star property.
The conjecture was proved independently by Ahlswede and
Khachatrian \cite{AK2} and Frankl and Tokushige \cite{FT2}
(Kleitman \cite{Kl} had long established this result for $k =
2$). As in Theorem~\ref{AK}, Ahlswede and Khachatrian \cite{AK2}
also determined the extremal structures.

\begin{theorem}[\cite{AK2}] \label{AKFT} Let $1 \leq t \leq n$ and
$k \geq 2$. Let $m$ be the largest integer such that $t + 2m <
\min\{n+1, t + 2\frac{t-1}{k-2}\}$ (by convention,
$\frac{t-1}{k-2} = \infty$ if $k=2$). \\
(i) If $(k,t) \neq (2,1)$ and $\frac{t-1}{k-2}$ is not integral,
then $\mathcal{A}$ is a largest $t$-intersecting sub-family of
$\mathcal{S}_{[n],k}$ if and only if
\[\mathcal{A} = \{A \in \mathcal{S}_{[n],k} \colon |A \cap X| \geq
t + m\}\]
for some $X \in \mathcal{S}_{Y,k}$ with $Y \in {[n] \choose
t+2m}$.\\
(ii) If $(k,t) \neq (2,1)$ and $\frac{t-1}{k-2}$ is integral, then
$\mathcal{A}$ is a largest $t$-intersecting sub-family of
$\mathcal{S}_{[n],k}$ if and only if
\[\mathcal{A} = \{A \in \mathcal{S}_{[n],k} \colon |A \cap X| \geq
t + j\}\]
for some $j \in \{m,m+1\}$ and some $X \in \mathcal{S}_{Y,k}$ with
$Y \in {[n] \choose t+2j}$.
\\
(iii) If $(k,t) = (2,1)$, then $\mathcal{A}$ is a largest
$t$-intersecting sub-family of $\mathcal{S}_{[n],k}$ if and only
if for any $y_1, ..., y_n \in [2]$, exactly one of $\{(1,y_1),
..., (n,y_n)\}$ and $\{(1,3-y_1), ..., (n,3-y_n)\}$ is in
$\mathcal{A}$.
\end{theorem}
Note that (iii) follows trivially from the fact that for any set
$A := \{(1,y_1), ..., (n,y_n)\}$ in $\mathcal{S}_{[n],2}$,
$\{(1,3-y_1), ..., (n,3-y_n)\}$ is the only set in
$\mathcal{S}_{[n],2}$ that does not intersect $A$. The rest of the
theorem is highly non-trivial!

%
What led to Theorem~\ref{AKFT} was the accomplishment of
Theorem~\ref{AK}. The following is an immediate consequence of
Theorem~\ref{AKFT}.

\begin{cor}\label{AKFTcor} Let $1 \leq t \leq n$ and
$k \geq 2$. Then:\\
(i) $\mathcal{S}_{[n],k}$ has the $t$-star property if and only if
$k \geq t+1$;\\
(ii) $\mathcal{S}_{[n],k}$ has the strict $t$-star property if and
only if $k \geq t+2$.
\end{cor}

We point out that Bey and Engel \cite{BE} extended
Theorem~\ref{AKFT} by determining the size of a largest
non-trivial $t$-intersecting sub-family of $\mathcal{S}_{[n],k}$
(see Examples 10, 11 and Lemma 18 in \cite{BE}).

Note that $\mathcal{S}_{[n],k} = \mathcal{S}_{{[n] \choose r},k}$
with $r = n$. For the case $t \leq r < n$, Bey \cite{Bey1} proved
the following.

\begin{theorem}[\cite{Bey1}] \label{Bey1} Let $1 \leq t \leq r
< n$. $\mathcal{S}_{{[n] \choose r},k}$ has the $t$-star property
if and only if $n \geq \frac{(r - t + k)(t+1)}{k}$.
\end{theorem}
Thus, if $t \leq r < n$ and $n \geq \frac{(r - t + k)(t+1)}{k}$,
then the size of a $t$-intersecting sub-family of
$\mathcal{S}_{{[n] \choose r},k}$ is at most ${n-t \choose
r-t}k^{r-t}$, i.e.~the size of any $t$-star of $\mathcal{S}_{{[n]
\choose r},k}$. From Corollary~\ref{AKFTcor} and
Theorem~\ref{Bey1} we immediately obtain the following.

\begin{cor} \label{Beycor} For any $1 \leq t \leq r \leq n$ and $k
\geq t+1$, $\mathcal{S}_{{[n] \choose r},k}$ has the $t$-star
property.
\end{cor}

To the best of the author's knowledge, no complete
$t$-intersection theorem for $\mathcal{S}_{{[n] \choose r},k}$
has been obtained.

For the case when $\mathcal{F}$ is any family, the author
\cite{Borg2} suggested the following general conjecture.

\begin{conj}[\cite{Borg2}] \label{ss_conj_t} For any integer $t
\geq 1$, there exists an integer $k_0(t)$ such that for any $k
\geq k_0(t)$ and any family $\mathcal{F}$,
$\mathcal{S}_{\mathcal{F},k}$ has the $t$-star property.
\end{conj}
In view of Corollary~\ref{Beycor}, we conjecture that the
smallest $k_0(t)$ is $t+1$. In \cite{Borg2} it is actually
conjectured that for some integer $k_0'(t)$,
$\mathcal{S}_{\mathcal{F},k}$ has the strict $t$-star property
for any $\mathcal{F}$, and hence, in view of
Corollary~\ref{AKFTcor}(ii), we conjecture that the smallest
$k_0'(t)$ is $t+2$. Note that the conjectures we have made about
the smallest values of $k_0(t)$ and $k_0'(t)$ generalise
Conjecture~\ref{ss conj}. The author \cite{Borg2} proved the
following relaxation of the statement of
Conjecture~\ref{ss_conj_t}.

\begin{theorem}[\cite{Borg2}] \label{ss_t_int2} For any integers
$r$ and $t$ with $1 \leq t < r$, let $k_0(r,t) := {r \choose t}{r
\choose t+1}$. For any $k \geq k_0(r,t)$ and any family
$\mathcal{F}$ with $\alpha(\mathcal{F}) \leq r$,
$\mathcal{S}_{\mathcal{F},k}$ has the strict $t$-star property.
\end{theorem}
The general idea behind the proof of this result is similar to
that behind the proof of Theorem~\ref{t int her}, described in
Section~\ref{Hersection}.

\begin{cor} Conjecture \ref{ss conj} is true if $k \geq
\alpha(\mathcal{F}){\alpha(\mathcal{F}) \choose 2}$.
\end{cor}

\section{Intersecting families of labeled sets} \label{Label_Int}

Consider the family $\mathcal{L}_{{\bf k}}$, ${\bf k} = (k_1,
..., k_n)$, of \emph{labeled $n$-sets}. If $k_1 = 1$ then all the
sets in $\mathcal{L}_{{\bf k}}$ contain the point $(1,1)$ and
hence $\mathcal{L}_{{\bf k}}$ has the strict star property. Berge
\cite{B} proved that for any ${\bf k}$, $\mathcal{L}_{{\bf k}}$
has the star property, and hence the size of an intersecting
sub-family of $\mathcal{L}_{{\bf k}}$ is at most the size
$\frac{1}{k_1}|\mathcal{L}_{{\bf k}}| = k_2k_3...k_n$ of the star
$\mathcal{L}_{{\bf k}}((1,1))$, as this is clearly a largest star
(since $k_1 \leq ... \leq k_n$). We shall reproduce the remarkably
short proof of this result.

Let $\mbox{mod}^*$ be the usual modulo operation with the
exception that for any integer $a$, $a \; \mbox{mod}^* \; a$ is
$a$ instead of $0$. For any integer $q$, let $\theta_{\bf k}^q :
\mathcal{L}_{\bf{k}} \rightarrow \mathcal{L}_{\bf{k}}$ be the
\emph{translation operation} defined by
\begin{equation}\theta_{\bf k}^q(A) := \{(a, (b+q) \; \mbox{mod}^*
\; k_a) \colon (a,b) \in A\}, \nonumber
\end{equation}
and define $\Theta_{\bf k}^q : 2^{\mathcal{L}_{\bf{k}}}
\rightarrow 2^{\mathcal{L}_{\bf{k}}}$ by
\begin{equation}\Theta_{\bf k}^q(\mathcal{F}) := \{\theta_{\bf
k}^q(A) \colon A \in \mathcal{F}\}.  \nonumber
\end{equation}
Let $\mathcal{A}$ be an intersecting sub-family of
$\mathcal{L}_{\bf k}$. For any $A \in \mathcal{A}$ and $q \in [k_1
- 1]$, we have $\theta_{\bf k}^q(A) \cap A = \emptyset$ and hence
$\theta_{\bf k}^q(A) \notin \mathcal{A}$. Therefore $\mathcal{A},
\Theta_{\bf k}^1(\mathcal{A}), ..., \Theta_{\bf
k}^{k_1-1}(\mathcal{A})$ are $k_1$ disjoint sub-families of
$\mathcal{L}_{\bf k}$. So $k_1|\mathcal{A}| \leq |\mathcal{L}_{\bf
k}|$ and hence $|\mathcal{A}| \leq \frac{1}{k_1}|\mathcal{L}_{\bf
k}|$.

Livingston \cite{L} proved that for $3 \leq k_1 = ... = k_n$,
$\mathcal{L}_{\bf k}$ has the strict star property. Using the
shifting technique (see \cite{F}) in an inductive argument, the
author \cite{Borg4} extended Livingston's result for the case when
$3 \leq k_1 \leq ... \leq k_n$. The above results sum up as
follows.

\begin{theorem}[\cite{B,Borg4,L}]\label{Liv_ext} Let $1 \leq k_1
\leq ... \leq k_n$ and let ${\bf k} := (k_1, ..., k_n)$. Then: \\
(i) $\mathcal{L}_{\bf{k}}$ has the star property at $\{(1,1)\}$;\\
(ii) if $k_1 \neq 2$ then $\mathcal{L}_{\bf{k}}$ has the strict
star property.
\end{theorem}
If $k_1 = 2$ then $\mathcal{L}_{\bf{k}}$ may not have the strict
star property; indeed, if $k_1 = k_2 = k_3$ then $\{A \in
\mathcal{L}_{\bf{k}} \colon |A \cap \{(1,1),(2,1),(3,1)| \geq
2\}$ is a non-trivial intersecting sub-family of
$\mathcal{L}_{\bf k}$ whose size is
$\frac{1}{k_1}|\mathcal{L}_{\bf k}|$ (i.e.~the maximum).

Recall that $\mathcal{S}_{[n],k} = \mathcal{L}_{(k_1, ..., k_n)}$
with $k_1 = ... = k_n = k$. The same argument used in \cite{Borg4}
to extend Livingston's result \cite{L} gives the following
extension of part (the sufficiency conditions) of
Corollary~\ref{AKFTcor} and generalisation of
Theorem~\ref{Liv_ext} with $k_1 \geq 2$.

\begin{theorem}\label{Liv_ext2} Let $2 \leq t+1
\leq k_1 \leq ... \leq k_n$ and let ${\bf k} := (k_1, ..., k_n)$.
Then: \\
(i) $\mathcal{L}_{\bf{k}}$ has the $t$-star property at $\{(1,1),
..., (t,1)\}$;\\
(ii) if $k_1 \geq t+2$ then $\mathcal{L}_{\bf{k}}$ has the strict
$t$-star property.
\end{theorem}
As we can see from Theorem~\ref{AKFT} and
Corollary~\ref{AKFTcor}, $\mathcal{L}_{\bf{k}}$ may not have the
$t$-star property when $2 \leq k_1 \leq t$. Recall that for the
case $k_1 = ... = k_n$, the extremal structures are given in
Theorem~\ref{AKFT}, and they are all non-trivial when $2 \leq k_1
\leq t$.

The intersection problem for the families $\mathcal{L}_{{\bf
k},r}$, $r = 1, ..., n$, has also been treated to a significant
extent. Note that $\mathcal{S}_{{[n] \choose r},k} =
\mathcal{L}_{(k_1, ..., k_n),r}$ with $k_1 = ... = k_n = k$.
Using the shifting technique (see \cite{F}) in an inductive
argument, Holroyd, Spencer and Talbot \cite{HST} extended
Theorem~\ref{DFBL}(i) as follows.

\begin{theorem}[\cite{HST}]\label{HSTthm} Let $2 \leq k_1 \leq ...
\leq k_n$ and let ${\bf k} := (k_1, ..., k_n)$. Then for any $r
\in [n]$, $\mathcal{L}_{{\bf k},r}$ has the star property at
$\{(1,1)\}$.
\end{theorem}
The proof of their result can be easily extended to obtain that
$\mathcal{L}_{{\bf k},r}$ has the strict star property if $(r,k_1)
\neq (n,2)$ (see, for example, the proof of
\cite[Theorem~1.4]{Borg4}). The case $k_1 = 1$ proved to be
harder, and Bey \cite{Bey2} solved it by applying the idea of
\emph{generating sets} introduced in \cite{AK1}.

\begin{theorem}[\cite{Bey2}]\label{Beyweight} Let $1 = k_1 = ...
= k_m < k_{m+1} \leq ... \leq k_n$ and let ${\bf k} := (k_1, ...,
k_n)$. Let $p := \lfloor (m+1)/2 \rfloor$, and for each $i \in
[p]$, let $\mathcal{A}_i := \{A \in \mathcal{L}_{{\bf k},r}
\colon (1,1) \in A, \, i \leq |A \cap \{(1,1), ..., (m,1)\}| \leq
m-i\} \cup \{A \in \mathcal{L}_{{\bf k},r} \colon |A \cap
\{(1,1), ..., (m,1)\}| \geq m-i+1\}$. Then one of the families
$\mathcal{A}_1, ..., \mathcal{A}_p$ is a largest intersecting
sub-family of $\mathcal{L}_{{\bf k},r}$.
\end{theorem}
Bey \cite{Bey2} also showed that when $r \leq n/2$ in the above
theorem, $\mathcal{L}_{{\bf k},r}$ has the star property at
$(1,1)$ (this is also proved in \cite{HST}, and in \cite{BH} it
is shown that $\mathcal{L}_{{\bf k},r}$ has the strict star
property if $r < n/2$).

For the case when $k_1$ can be any positive integer but $n$ is
sufficiently large, Theorem~\ref{t int her} gives us the
following $t$-intersection result.

\begin{theorem}\label{t-intlabel} Let $1 \leq t \leq r$ and let
$n \geq (r-t){3r-2t-1 \choose t+1} + r$. Let $1 \leq k_1 \leq ...
\leq k_n$ and let ${\bf k} := (k_1, ..., k_n)$. Then:\\
(i) $\mathcal{L}_{{\bf k},r}$ has the $t$-star property at
$\{(1,1), ..., (t,1)\}$.\\
(ii) $\mathcal{L}_{{\bf k},r}$ has the strict $t$-star
property.
\end{theorem}
\textbf{Proof.} Let $\mathcal{H} = \mathcal{L}_{{\bf k}, \leq
n}$. Then clearly $\mathcal{H}$ is a hereditary family with
$\mu(\mathcal{H}) = n$. Thus, by Theorem~\ref{t int her} (with $S
= \{r\}$), $\mathcal{H}^{(r)}$ has the strict $t$-star property.
Part (ii) follows since $\mathcal{H}^{(r)} = \mathcal{L}_{{\bf
k},r}$. This in turn proves (i) since the family
$\mathcal{L}_{{\bf k},r}(T)$ with $T := \{(1,1), ..., (t,1)\}$ is
clearly a largest $t$-star of
$\mathcal{L}_{{\bf k},r}$.~\hfill{$\Box$}\\

We mention that Erd\H os, Seress, and Sz\'{e}kely \cite{EES}
determined non-trivial $t$-intersecting sub-families of
$\mathcal{L}_{{\bf k},r}$ of maximum size for the case when $n$
is sufficiently large.

Finally, for the family $\mathcal{L}_{{\bf k},\leq n}$ of all
labeled sets defined on the $n$-tuple ${\bf k}$, we have the
following immediate consequence of Theorems~\ref{Snevily}
and~\ref{HSTthm}.

\begin{theorem}\label{leqn} For any $1 \leq k_1 \leq ... \leq
k_n$, $\mathcal{L}_{(k_1, ..., k_n),\leq n}$ has the star property
at $\{(1,1)\}$.
\end{theorem}
\textbf{Proof.} Let ${\bf k} = (k_1, ..., k_n)$. If $k_1 = 1$
then $\mathcal{L}_{{\bf k},\leq n}$ is compressed with respect to
$(1,1)$ and hence, since $\mathcal{L}_{{\bf k},\leq n}$ is
hereditary, the result follows by Theorem~\ref{Snevily}. Now
suppose $k_1 \geq 2$. Let $\mathcal{A}$ be an intersecting
sub-family of $\mathcal{L}_{{\bf k},\leq n}$. So $\emptyset
\notin \mathcal{A}$. By Theorem~\ref{HSTthm}, $|\mathcal{A}^{(r)}|
\leq |\mathcal{L}_{{\bf k},r}((1,1))|$ for all $r \in [n]$. Thus,
we have $|\mathcal{A}| = \sum_{r=1}^n |\mathcal{A}^{(r)}| \leq
\sum_{r=1}^n |\mathcal{L}_{{\bf k},r}((1,1))| =
|\mathcal{L}_{{\bf k},\leq n}((1,1))|$.~\hfill{$\Box$}\\

The above fact was also observed in \cite{Bey2}, and it implies
that the size of an intersecting sub-family of $\mathcal{L}_{{\bf
k},\leq n}$ is at most $\frac{1}{k+1}|\mathcal{L}_{{\bf k},\leq
n}|$, i.e.~the size of the star $\mathcal{L}_{{\bf k},\leq
n}((1,1))$ (indeed, the $k_1 + 1$ families $\mathcal{L}_{{\bf
k},\leq n}((1,1)), ..., \mathcal{L}_{{\bf k},\leq n}((1,k_1))$ and
$\mathcal{L}_{(k_2, ..., k_n),\leq n-1}$ partition
$\mathcal{L}_{{\bf k},\leq n}$ and are of the same size). In view
of the above-mentioned fact that $\mathcal{L}_{{\bf k},r}$ has
the strict star property when $(r,k_1) \neq (n,2)$ (in
particular, when $1 \leq r \leq n-1$ and $k_1 \geq 2$), one can go
on to show that $\mathcal{L}_{{\bf k},\leq n}$ has the strict star
property if $k_1 \geq 2$. If $k_1 = 1$ then $\mathcal{L}_{{\bf
k},\leq n}$ may not have the strict star property; indeed, if
$k_1 = k_2 = k_3 = 1$ then $\{A \in \mathcal{L}_{{\bf k},\leq n}
\colon |A \cap \{(1,1),(2,1),(3,1)\}| \geq 2\}$ is a non-trivial
intersecting sub-family that is as large as the largest star
$\mathcal{L}_{{\bf k},\leq n}((1,1))$.

To the best of the author's knowledge, no general $t$-intersection
theorem for $\mathcal{L}_{{\bf k}, \leq n}$ is known.

\section{Intersecting families of permutations and partial
permutations}

In \cite{Deza,DF1} the study of intersecting permutations was
initiated. Deza and Frankl \cite{DF1} showed that
$\mathcal{S}_{[n],n}^*$ has the star property. So the size of an
intersecting sub-family of $\mathcal{S}_{[n],n}^*$ is at most
$(n-1)!$. The argument of the proof of this result is the same
translation argument, given in the previous section, that yields
Berge's intersection result for labeled sets \cite{B}, and it
also gives us that for $n \leq k$, $\mathcal{S}_{[n],k}^*$ has
the star property (recall that $\mathcal{S}_{[n],k}^* =
\emptyset$ if $n > k$). Indeed, it gives us that for any
intersecting sub-family $\mathcal{A}$ of $\mathcal{S}_{[n],k}^*$,
$k|\mathcal{A}| \leq |\mathcal{S}_{[n],k}^*| = \frac{k!}{(k-n)!}$
and hence $|\mathcal{A}| \leq \frac{(k-1)!}{(k-n)!}$.

The question of whether $\mathcal{S}_{[n],n}^*$ has the strict
star property proved to be much more difficult to answer. Cameron
and Ku \cite{CK} and Larose and Malvenuto \cite{LM} independently
gave an affirmative answer (other proofs are given in
\cite{GM,WZ}). Larose and Malvenuto \cite{LM} also proved the
following generalisation (another proof is found in \cite{BrH}).

\begin{theorem}[\cite{LM}] For $1 \leq n \leq k$,
$\mathcal{S}_{[n],k}^*$ has the strict star property.
\end{theorem}

Ku and Leader \cite{KL} investigated partial permutations. Using
Katona's cycle method \cite{K}, they proved that
$\mathcal{S}_{{[n] \choose r},n}^*$ has the star property for all
$r \in [n-1]$ (note that if $r = n$, then $\mathcal{S}_{{[n]
\choose r},n}^* = \mathcal{S}_{[n],n}^*$), and they also showed
that $\mathcal{S}_{{[n] \choose r},n}^*$ has the strict star
property for all $r \in [8,n-3]$. Naturally, they conjectured
that $\mathcal{S}_{{[n] \choose r},n}^*$ also has the strict star
property for the few remaining values of $r$. This was settled by
Li and Wang \cite{LW} using tools forged by Ku and Leader. So the
intersection results for $\mathcal{S}_{[n],n}^*$ and
$\mathcal{S}_{{[n] \choose r},n}^*$ ($r \in [n-1]$) sum up as
follows.

\begin{theorem}[\cite{CK,KL,LM,LW}] For any $r \in [n]$,
$\mathcal{S}_{{[n] \choose r},n}^*$ has the strict star property.
\end{theorem}

When it comes to $t$-intersecting families of permutations, things
are of course much harder. Solving a long-standing conjecture of
Deza and Frankl \cite{DF1}, Ellis, Friedgut and Pilpel \cite{EFP}
recently managed to prove the following.

\begin{theorem}[\cite{EFP}]\label{DF_perm} For any integer
$t \geq 1$, there exists an integer $n_0(t)$ such
that for any $n \geq n_0(t)$, $\mathcal{S}_{[n],n}^*$ has the strict $t$-star property.
\end{theorem}
Their remarkable proof is based on eigenvalues techniques and
representation theory of the symmetric group. The condition $n
\geq n_0(t)$ is necessary. Indeed, let $P_j := \{(i,i) \colon i
\in [j]\}$ for any integer $j \geq 1$, and let
\[\mathcal{G}_{n,k,t} := \left\{ \begin{array}{ll}
\{A \in \mathcal{S}_{[n],k} \colon |A \cap P_n| \geq (n+t)/2\} & \mbox{if $n-t$ is even};\\
\{A \in \mathcal{S}_{[n],k} \colon |A \cap P_{n-1}| \geq
(n+t-1)/2\} & \mbox{if $n-t$ is odd}.
\end{array} \right. \]
Deza and Frankl \cite{DF1} showed that when $t = n-s$ for some $s
\geq 3$ and $n$ is sufficiently large, $\mathcal{G}_{n,k,t}$ is a
largest $t$-intersecting sub-family of $\mathcal{S}_{[n],n}^*$ and
is larger than the $t$-stars. Brunk and Huczynska \cite{BrH}
extended this result as follows.

\begin{theorem}[\cite{BrH,DF1}] For any integers $p \geq 2$ and $q
\geq 0$, there exists an integer $n_0^*(p,q)$ such that for any
$n \geq n_0^*(p,q)$, \\
(i) $\mathcal{G}_{n,n+p,n-q}$ is a largest $(n-q)$-intersecting
sub-family of $\mathcal{S}_{[n],n+p}^*$;\\
(ii) any largest $(n-q)$-intersecting sub-family of
$\mathcal{S}_{[n],n+p}^*$ is a copy of $\mathcal{G}_{n,n+p,n-q}$.
\end{theorem}
They also conjectured that for any $n$ and $k \geq n$, the
extremal structures are similar to those in Theorem~\ref{AK}.
\begin{conj}[\cite{BrH}] Let $1 \leq t \leq n \leq k$. Let $p :=
\lfloor (n-t)/2 \rfloor$, and for any integer $i$ with $0 \leq i
\leq p$, let
\[\mathcal{A}_i := \left\{A \in \mathcal{S}_{[n],k}^* \colon |A \cap
\{(1,1), ..., (t+2i,t+2i)\}| \geq t+i \right\}.\]
Then: \\
(i) one of the families $\mathcal{A}_0, ..., \mathcal{A}_p$ is a
largest $t$-intersecting sub-family of
$\mathcal{S}_{[n],k}^*$;\\
(ii) any largest $t$-intersecting sub-family of
$\mathcal{S}_{[n],k}^*$ is a copy of one of the families
$\mathcal{A}_0, ..., \mathcal{A}_p$.
\end{conj}

For the general case when $\mathcal{F}$ is any family, a
conjecture for $t$-intersecting sub-families of
$\mathcal{S}_{\mathcal{F},k}^*$ similar to
Conjecture~\ref{ss_conj_t} was suggested in \cite{Borg2}.
\begin{conj}[\cite{Borg2}] For any integer $t \geq 1$, there
exists an integer $k_0^*(t)$ such that for any $k \geq k_0^*(t)$
and any family $\mathcal{F}$, $\mathcal{S}_{\mathcal{F},k}^*$ has
the strict $t$-star property.
\end{conj}
Theorem~\ref{DF_perm} solves the special case $\mathcal{F} =
\{[n]\}$ and $k = n \geq k_0^*(t)$. The author \cite{Borg2}
proved the following relaxation of the statement of the
conjecture.

\begin{theorem}[\cite{Borg2}]\label{t_int_perm} For any integers $r$ and $t$
with $1 \leq t < r$, let $k_0^*(r,t) := {r \choose t}{3r-2t-1
\choose \lfloor \frac{3r-2t-1}{2} \rfloor} \frac{r!}{(r-t-1)!} +
r + 1$. For any $k \geq k_0^*(r,t)$ and any family $\mathcal{F}$
with $\alpha(\mathcal{F}) \leq r$,
$\mathcal{S}_{\mathcal{F},k}^*$ has the strict $t$-star property.
\end{theorem}
This is an analogue of Theorem~\ref{ss_t_int2}, and the general
idea behind its proof is similar to that behind the proofs of
Theorems~\ref{t int her} (see Section~\ref{Hersection})
and~\ref{ss_t_int2}.

By taking $\mathcal{F} = [n]$ and $k \geq k_0^*(n,t)$ in
Theorem~\ref{t_int_perm}, we obtain the following.

\begin{cor} Let $k \geq k_0^*(n,t)$, where $k_0^*(n,t)$ is as in Theorem~\ref{t_int_perm}. Then $\mathcal{S}_{[n],k}^*$ has the strict $t$-star property.
\end{cor}
Thus, when $k$ is sufficiently large, the size of a $t$-intersecting sub-family of $\mathcal{S}_{[n],k}^*$ is at most $\frac{(k-t)!}{(k-n)!}$.

The following $t$-intersection result for partial permutations is
another immediate consequence of Theorem~\ref{t_int_perm},
obtained by taking $n \geq k_0^*(r,t)$ and $\mathcal{F} = {[n]
\choose r}$.

\begin{cor} \label{Ku} Let $n \geq k_0^*(r,t)$, where $k_0^*(r,t)$ is as in Theorem~\ref{t_int_perm}. Then $\mathcal{S}_{{[n] \choose r},n}^*$ has the strict $t$-star property.
\end{cor}
Thus, when $n$ is sufficiently large, the size of a
$t$-intersecting sub-family of $\mathcal{S}_{{[n] \choose
r},n}^*$ is at most ${n-t \choose r-t}\frac{(n-t)!}{(n-r)!}$.
This was also proved in \cite{Ku_thesis}.


\begin{thebibliography}{}
\bibitem{AK1} R. Ahlswede, L.H. Khachatrian, The complete
intersection theorem for systems of finite sets, European J.
Combin. 18 (1997) 125-136.
\bibitem{AK3} R. Ahlswede, L.H. Khachatrian, The complete
nontrivial-intersection theorem for systems of finite sets, J.
Combin. Theory Ser. A 76 (1996) 121-138.
\bibitem{AK2} R. Ahlswede, L.H. Khachatrian, The diametric theorem
in Hamming spaces - Optimal anticodes, Adv. in Appl. Math. 20
(1998) 429-449.
\bibitem{Anderson} I. Anderson, Combinatorics of Finite Sets,
Oxford University Press, Oxford, England, 1987.
\bibitem{Berge} C. Berge, A theorem related to the Chvátal
conjecture, Proceedings of the Fifth British Combinatorial
Conference (Univ. Aberdeen, Aberdeen, 1975),  pp. 35-40.
Congressus Numerantium, No. XV, Utilitas Math., Winnipeg, Man.,
1976.
\bibitem{B} C. Berge, Nombres de coloration de l'hypergraphe
h-parti complet, in: Hypergraph Seminar (Columbus, Ohio 1972),
Lecture Notes in Math., Vol. 411, Springer, Berlin, 1974, pp.
13-20.
\bibitem{Bey2} C. Bey, An intersection theorem for weighted sets,
Discrete Math. 235 (2001) 145-150.
\bibitem{Bey1} C. Bey, The Erd\H{o}s-Ko-Rado bound for the function
lattice, Discrete Appl. Math. 95 (1999) 115-125.
\bibitem{BE} C. Bey, K. Engel, Old and new results for the weighted
t-intersection problem via AK-methods, In: Numbers, Information
and Complexity, Alth\H{o}fer, Ingo, Eds. et al., Dordrecht: Kluwer
Academic Publishers, 2000, pp. 45-74.
\bibitem{BL} B. Bollob\'{a}s, I. Leader, An Erd\H os-Ko-Rado theorem
for signed sets, Comput. Math. Appl. 34 (1997) 9-13.
\bibitem{Borg} P. Borg, Extremal $t$-intersecting sub-families of
hereditary families, J. London Math. Soc. 79 (2009) 167-185.
\bibitem{Borg4} P. Borg, Intersecting and cross-intersecting
families of labeled sets, Electron. J. Combin. 15 (2008), N9.
\bibitem{Borg1} P. Borg, Intersecting systems of signed sets,
Electron. J. Combin. 14 (2007) $\#$R41.
\bibitem{Borg3} P. Borg, On Chv\'{a}tal's conjecture and a
conjecture on families of signed sets, European J. Combin. 32
(2011) 140-145.
\bibitem{Borg2} P. Borg, On $t$-intersecting families of signed
sets and permutations, Discrete Math. 309 (2009) 3310-3317.
\bibitem{BH} P. Borg, F.C. Holroyd, The Erdos-Ko-Rado property of
various graphs containing singletons, Discrete Math. 309 (2009)
2877-2885.
\bibitem{BrH} F. Brunk, S. Huczynska, Some Erd\H{o}s-Ko-Rado
theorems for injections, European J. Combin. 31 (2010), 839-860.
\bibitem{CK} P.J. Cameron, C.Y. Ku, Intersecting families of
permutations, European J. Combin. 24 (2003) 881-890.
\bibitem{Chv} V. Chv\'{a}tal, Unsolved Problem No. 7, in: C. Berge,
D.K. Ray-Chaudhuri (Eds.), Hypergraph Seminar, Lecture Notes in
Mathematics, Vol. 411, Springer, Berlin, 1974.
\bibitem{Chva} V. Chv\'{a}tal, Intersecting families of edges in
hypergraphs having the hereditary property, in: C. Berge, D.K.
Ray-Chaudhuri (Eds.), Hypergraph Seminar, Lecture Notes in
Mathematics, Vol. 411, Springer, Berlin, 1974, pp. 61-66.
\bibitem{Chvatalsite} V. Chv\'{a}tal,
http://users.encs.concordia.ca/$\sim$chvatal/conjecture.html.
\bibitem{D} D.E. Daykin, Erd\H os-Ko-Rado from Kruskal-Katona,
J. Combin. Theory Ser. A 17 (1974) 254-255.
\bibitem{Deza} M. Deza, Matrices dont deux lignes quelconques
coincident dans un nombre donne de positions communes, J. Combin.
Theory Ser. A 20 (1976) 306-318.
\bibitem{DF1} M. Deza, P. Frankl, On the maximum number of
permutations with given maximal or minimal distance, J. Combin.
Theory Ser. A 22 (1977) 352-360.
\bibitem{DF} M. Deza, P. Frankl, The Erd\H os-Ko-Rado theorem - 22
years later, SIAM J. Algebraic Discrete Methods 4 (1983) 419-431.
\bibitem{EFP} D. Ellis, E. Friedgut, H. Pilpel, Intersecting families of permutations, J. Amer. Math. Soc., to appear.
\bibitem{E} K. Engel, An Erd\H os-Ko-Rado theorem for the subcubes
of a cube, Combinatorica 4 (1984) 133-140.
\bibitem{EKR} P. Erd\H os, C. Ko and R. Rado, Intersection
theorems for systems of finite sets, Quart. J. Math. Oxford (2)
12 (1961) 313-320.
\bibitem{EFK} P.L. Erd\H os, U. Faigle, W. Kern, A group-theoretic
setting for some intersecting Sperner families, Combin. Probab.
Comput. 1 (1992) 323-334.
\bibitem{EES} P.L. Erd\H os, A. Seress, L.A.
Sz\'{e}kely, Non-trivial $t$-intersection in the function
lattice, Ann. Combin. 9 (2005) 177-187.
\bibitem{F_t1} P. Frankl, The Erd\H os-Ko-Rado Theorem is true
for $n = ckt$, Proc. Fifth Hung. Comb. Coll., North-Holland,
Amsterdam, 1978, pp. 365-375.
\bibitem{F} P. Frankl, The shifting technique in extremal set
theory, in: C. Whitehead (Ed.), Combinatorial Surveys, Cambridge
Univ. Press, London/New York, 1987, pp. 81-110.
\bibitem{FF2} P. Frankl, Z. Füredi, The Erd\H os-Ko-Rado Theorem
for integer sequences, SIAM J. Algebraic Discrete Methods 1(4)
(1980) 376-381.
\bibitem{FT2} P. Frankl, N. Tokushige, The Erd\H os-Ko-Rado
theorem for integer sequences, Combinatorica 19 (1999) 55-63.
\bibitem{GM} C. Godsil, K. Meagher, A new proof of the
Erd\H{o}s-Ko-Rado theorem for intersecting families of
permutations, European J. Combin. 30 (2009) 404-414.
\bibitem{G} H.-D.O.F. Gronau, More on the Erd\H os-Ko-Rado theorem
for integer sequences, J. Combin. Theory Ser. A 35 (1983) 279-288.
\bibitem{HM} A.J.W. Hilton, E.C. Milner, Some intersection
theorems for systems of finite sets, Quart. J. Math. Oxford (2)
18 (1967) 369-384.
\bibitem{HST} F.C. Holroyd, C. Spencer, J. Talbot, Compression and
Erd\H os-Ko-Rado graphs, Discrete Math. 293 (2005) 155-164.
\bibitem{HT} F.C. Holroyd, J. Talbot, Graphs with the Erd\H
os-Ko-Rado property, Discrete Math. 293 (2005) 165-176.
\bibitem{K} G.O.H. Katona, A simple proof of the Erd\H os-Chao
Ko-Rado theorem, J. Combin. Theory Ser. B 13 (1972) 183-184.
\bibitem{Ka} G.O.H. Katona, A theorem of finite sets, in:
Theory of Graphs, Proc. Colloq. Tihany, Akadémiai Kiadó (1968)
187-207.
\bibitem{Kat} G.O.H. Katona, Intersection theorems for finite sets,
Acta Math. Acad. Sci. Hungar. 15 (1964) 329-337.
\bibitem{Kl} D.J. Kleitman, On a combinatorial conjecture of Erd\H
os, J. Combin. Theory Ser. A 1 (1966) 209-214.
\bibitem{Kr} J.B. Kruskal, The number of simplices in a
complex, in: Mathematical Optimization Techniques, University of
California Press, Berkeley, California, 1963, pp. 251-278.
\bibitem{Ku_thesis} C.Y. Ku, Intersecting families of permutations
and partial permutations, Ph.D. Dissertation, Queen Mary College,
University of London, December, 2004.
\bibitem{KL} C.Y. Ku, I. Leader, An Erd\H os-Ko-Rado theorem for
partial permutations, Discrete Math. 306 (2006) 74-86.
\bibitem{LM} B. Larose, C. Malvenuto, Stable sets of maximal size
in Kneser-type graphs, European J. Combin. 25 (2004) 657-673.
\bibitem{LW} Y.-S. Li, Jun Wang, Erd\H os-Ko-Rado-type
theorems for colored sets, Electron. J. Combin. 14 (2007) $\#$R1.
\bibitem{L} M.L. Livingston, An ordered version of the Erd\H
os-Ko-Rado Theorem, J. Combin. Theory Ser. A 26 (1979), 162-165.
\bibitem{Meyer} J.-C. Meyer, Quelques problèmes concernant les
cliques des hypergraphes $k$-complets et $q$-parti $h$-complets,
in: Hypergraph Seminar (Columbus, Ohio 1972), Lecture Notes in
Math., Vol. 411, Springer, Berlin, 1974, 127-139.
\bibitem{Miklos} D. Mikl\'{o}s, Some results related to a
conjecture of Chv\'{a}tal, Ph.D. Dissertation, Ohio State
University, 1986.
\bibitem{M} A. Moon, An analogue of the Erd\H os-Ko-Rado theorem
for the Hamming schemes $H(n,q)$, J. Combin. Theory Ser. A 32
(1982) 386-390.
\bibitem{Schonheim} J. Schönheim, Hereditary systems and Chv\'{a}tal's
conjecture, Proceedings of the Fifth British Combinatorial
Conference (Univ. Aberdeen, Aberdeen, 1975), pp. 537-539.
Congressus Numerantium, No. XV, Utilitas Math., Winnipeg, Man.,
1976.
\bibitem{Sn} H.S. Snevily, A new result on Chv\'{a}tal's conjecture, J.
Combin. Theory Ser. A 61 (1992) 137-141.
\bibitem{Wang} D.L. Wang and P. Wang, Some results about the
Chv\'{a}tal conjecture, Discrete Math. 24 (1978) 95-101.
\bibitem{WZ} J. Wang, S.J. Zhang, An Erd\H{o}s-Ko-Rado-type
theorem in Coxeter groups, European J. Combin. 29 (2008)
1112-1115.
\bibitem{West} D.B. West,
http://www.math.uiuc.edu/$\sim$west/regs/chvatal.html.
\bibitem{W} R.M. Wilson, The exact bound in the Erd\H os-Ko-Rado
theorem, Combinatorica 4 (1984) 247-257.
\end{thebibliography}
\end{document}